\documentclass[11pt]{amsart}
\usepackage[T1]{fontenc}
\usepackage{lmodern,microtype}
\usepackage[margin=1in]{geometry}
\usepackage{amssymb,mathtools,xcolor,needspace,tikz,placeins}
\usepackage[colorlinks=true,linkcolor=blue!45!black,citecolor=blue!45!black,urlcolor=blue!45!black]{hyperref}
\hypersetup{pdftitle={Real Hadamard matrices and Temperley--Lieb subfactors},pdfsubject={Hadamard matrices and relative commutants},pdflang={en}}
\newtheorem{theorem}{Theorem}[section]
\newtheorem{proposition}[theorem]{Proposition}
\newtheorem{lemma}[theorem]{Lemma}
\newtheorem{corollary}[theorem]{Corollary}
\theoremstyle{remark}

\numberwithin{equation}{section}
\newcommand{\C}{\mathbb C}
\newcommand{\Q}{\mathbb Q}

\newcommand{\Tr}{\operatorname{Tr}}
\newcommand{\Cat}{\operatorname{Cat}}
\newcommand{\TL}{\mathrm{TL}}
\newcommand{\cH}{\mathcal H}
\newcommand{\cP}{\mathcal P}
\newcommand{\cQ}{\mathcal Q}
\newcommand{\cA}{\mathcal A}
\title[Real Hadamard matrices and subfactors]{Subfactor Symmetries of real 
Hadamard matrices} 
\author{Dietmar Bisch}
\address[Dietmar Bisch]{Department of Mathematics, Vanderbilt University, Nashville}
\thanks{Dietmar Bisch was supported in part by U.S. Army Research Office Grant W911NF-23-1-0026.}
\author{Remus Nicoara}
\address[Remus Nicoara]{Department of Mathematics, University of Tennessee, Knoxville}
\date{}
\begin{document}

\begin{abstract}
Every Hadamard matrix gives rise to a {\it spin model commuting square}
and a hyperfinite subfactor with finite Jones index. The standard 
invariant of this subfactor is comprised of quantum symmetries 
of the subfactor which can be viewed as ``non-invertible'' symmetries
of the Hadamard matrix. We show that if the Hadamard matrix is real 
and has order greater than four and congruent to four modulo eight, 
the standard invariant is minimal, i.e. it consists of Temperley-Lieb
algebras. Consequently, both principal graphs of the subfactor are 
$A_\infty$. In particular, this holds for
every Paley type II Hadamard matrix. The proof uses a congruence for
products of four rows and Jones's annular decomposition. More generally, we show that if an irreducible subfactor of finite index greater
than four has bounded annular multiplicities along an infinite
subsequence of degrees, then its principal graph is $A_\infty$.
\end{abstract}
\maketitle

\section{Introduction}

A real Hadamard matrix of order $d$ is a matrix $H\in M_d(\C)$
whose entries are $1$ or $-1$ and whose rows are mutually orthogonal.
Thus $HH^{\mathsf T}=dI_d$, and $U=d^{-1/2}H$ is unitary.
Let $D_d\subset M_d(\C)$ be the algebra of diagonal matrices.
The inclusions
\[
\begin{matrix}
D_d&\subset&M_d(\C)\\
\cup&&\cup\\
\C&\subset&UD_dU^*
\end{matrix}
\]
form a commuting square, called a {\it spin model} commuting square, 
with respect to the normalized trace \cite[Section~2.11]{J99}. This 
means that the two diagonal subalgebras of $M_d(\C)$ in the above
quadrilateral are orthogonal modulo the constant matrices with respect
to the inner product on $M_d(\C)$ defined by the trace.
 
By iterating Jones' basic construction, one obtains an irreducible 
inclusion
$N\subset M$ of hyperfinite $\mathrm{II}_1$ factors with Jones index $d$
\cite{J83}.
We use the convention of \cite[Section~2]{N10} for this construction.
The opposite convention replaces $H$ by $H^{\mathsf T}$ and does
not affect the hypotheses considered here.

Write
\[
N=M_{-1}\subset M=M_0\subset M_1\subset M_2\subset\cdots
\]
for the Jones tower, and let $e_k\in M_k$, $k\geq1$, be the Jones
projections. The {\it standard invariant} of the subfactor
$N \subset M$ consists of the relative commutants
$N'\cap M_k$ and $M'\cap M_k$, together with the data describing
their inclusions and traces \cite{GHJ89,P94}. 
The relative commutants are finite dimensional
C$^*$-algebras of intertwiners of tensor powers of
the standard $N$-$M$ bimodule L$^2$(M) with its contragredient.
These tensor powers lead to unitary tensor categories of $N$-$N$ and 
$M$-$M$ bimodules that are viewed as {\it quantum symmetries}
of the subfactor $N \subset M$. Thus, these categories of subfactor 
representations constitute invariants of the Hadamard 
matrix that goes into the construction of the subfactor.

Each centralizer algebra $N'\cap M_k$ contains the unital algebra generated
by the Jones projections $e_1,\ldots,e_k$. For $d>4$, these algebras 
are Temperley--Lieb algebras with parameter $\sqrt d$ \cite{GHJ89,J83}. 
It is therefore natural to ask whether the relative commutants contain 
any additional elements. 

This turns out to be a notoriously difficult question. While the standard
invariant of any subfactor that arises from a commuting square can
in principle be computed by solving precise systems of linear equations,
due to a result known as Ocneanu compactness (see e.g. \cite[Section~5.7]{JS97}),
it is in practice nearly impossible to carry out the computations. This
is due to the fact that the coefficient matrices for the system of
equations are hard to compute, and the number of unknowns tends to 
increase exponentially.

Jones developed the theory of planar algebras as a tool to compute
with the quantum symmetries of a subfactor. He was particularly 
interested in computing them for Hadamard (or spin model) subfactors
and showed that the profile matrix can be used to determine the
dimension of $N' \cap M_1$
in \cite[Theorem~2.11.15]{J99}. We are using planar algebra techniques in this
paper to show that hyperfinite subfactors associated to a certain
family of real Hadamard matrices have {\it minimal}
standard invariant, i.e. the higher relative commutants are generated 
by the Jones projections and hence are Temperley-Lieb algebras. 
The principal graphs of these subfactors are therefore the Coxeter-Dynkin
graph $A_\infty$.

The problem of determining the higher relative commutants of
a spin model subfactor was studied in \cite{N10} by computing the first few
relative commutants for Hadamard matrices of small orders and by
considering special classes of matrices. For the Paley type II
matrices, Montgomery proved that the associated subfactors have
infinite depth and conjectured that their standard invariants are
minimal \cite[Example~3.4 and Conjecture~3.5]{M21}.
In this paper we prove this assertion for every real Hadamard matrix
of order $d>4$ satisfying $d\equiv4\pmod8$.

We formulate the result using the associated shaded planar algebra
$P$. Our indexing convention is
\[
P_{n,+}=N'\cap M_{n-1},\qquad
P_{n,-}=M'\cap M_n\quad(n\geq1),\qquad
P_{0,+}=P_{0,-}=\C.
\]
This planar algebra is connected, irreducible and spherical, with
modulus $\delta=\sqrt d$. We write $P_n=P_{n,+}$ and denote by
$\TL_{n,\pm}(\delta)$ the canonical Temperley--Lieb subspaces.
The $n$th Catalan number is denoted by
$\Cat_n=(n+1)^{-1}\binom{2n}{n}$.

\begin{theorem}\label{thm:main}
Let $H$ be a real Hadamard matrix of order $d>4$ such that
$d\equiv4\pmod8$. Then the planar algebra of the associated subfactor
is Temperley--Lieb:
\[
P_{n,\pm}=\TL_{n,\pm}(\sqrt d)\qquad(n\geq0).
\]
Its principal graphs are $A_\infty$, and
\[
\dim(N'\cap M_k)=\Cat_{k+1},\qquad
\dim(M'\cap M_k)=\Cat_k\qquad(k\geq0).
\]
\end{theorem}

\begin{corollary}\label{cor:paley}
The subfactor associated to a Paley type II Hadamard matrix has
Temperley--Lieb standard invariant.
\end{corollary}
\begin{proof}
A Paley type II matrix has order $d=2(q+1)$, where $q$ is a prime
power congruent to one modulo four; see \cite[Example~3.4]{M21}.
It follows that $d>4$ and $d\equiv4\pmod8$.
The result is therefore a consequence of Theorem~\ref{thm:main}.
\end{proof}

We briefly describe the main ingredients of the proof. The first step
is to simplify the entries of the profile matrix using the 
classical four-row congruence recorded by Wallis~\cite[p.~425]{W89}
and used by Montgomery in the proof of~\cite[Proposition~3.2.9]{M22}:
\[
\sum_y H_{ay}H_{by}H_{cy}H_{ey}\equiv-d\pmod8.
\]
When $d\equiv4\pmod8$, this congruence implies the two-click rotation acts trivially on the common kernel
of the cap maps in every odd degree. A trace estimate, together with
Jones's annular decomposition, then gives a uniform bound on the
odd-degree annular multiplicities. From this bound we deduce that
the principal graph has norm at most two. The classification of such
graphs and Popa's theorem on principal graphs with an $A_\infty$
tail~\cite{P95} then identify the graph as $A_\infty$.

The paper is organized as follows. In Section~\ref{sec:cyclic} we
express the relative commutants, rotation and cap maps in coordinates
coming from the finite dimensional algebras in the commuting-square
tower. Section~\ref{sec:parity} contains the congruence argument, and
Section~\ref{sec:trace} gives the trace estimate.
In Section~\ref{sec:annular} we recall the results of Jones on annular
Temperley--Lieb modules that we need. In the last section we prove
the graph-norm criterion and complete the proof of the main theorem.

\medskip
\noindent\textbf{AI Statement.}
OpenAI's ChatGPT assisted with the development and exposition of the
arguments and with checking references. The authors take responsibility
for the mathematical content of the paper.

\section{Cyclic coordinates for the relative commutants}\label{sec:cyclic}
Let $H$ be a real Hadamard matrix of order $d$, with rows and columns
indexed by a set $X$. No congruence assumption on $d$ is needed in this
section. For $a,b,c,e\in X$, put
\begin{equation}\label{eq:profile}
 S(a,b,c,e)=\sum_{y\in X}H_{ay}H_{by}H_{cy}H_{ey}.
\end{equation}
These sums determine the profile matrix and the higher relative
commutants. To describe the latter, let $\cH_n=\ell^2(X^n)$, with
inner product
\[
\langle f,g\rangle=\sum_{a\in X^n}f(a)\overline{g(a)}.
\]
We read subscripts cyclically, so that $a_{n+1}=a_1$.
For $n\geq1$, define matrices $C_n,T_n$ on $\cH_n$ by
\begin{align}
 C_n(a,b)&=d^{-(n+1)/2}\prod_{j=1}^n H_{a_jb_j}H_{a_{j+1}b_j},\label{eq:C}\\
 T_n(a,b)&=d^{-n-1}\prod_{j=1}^n S(a_j,b_j,a_{j+1},b_{j+1}),\label{eq:T}
\end{align}
and let
\begin{equation}\label{eq:F}
 F_n=\ker(I-T_n),\qquad
 \rho e_{(a_1,\ldots,a_n)}=e_{(a_2,\ldots,a_n,a_1)}.
\end{equation}
Equivalently, on coefficient functions,
\[
 (\rho f)(a_1,\ldots,a_n)=f(a_n,a_1,\ldots,a_{n-1}),
 \qquad \rho^n=I.
\]
For $n\ge2$, introduce the orthogonal projections
\begin{align*}
 (A_jf)(a)&=\frac1d\sum_{x\in X}f(a_1,\ldots,a_{j-1},x,a_{j+1},\ldots,a_n),\\
 (E_jf)(a)&=\mathbf1_{a_j=a_{j+1}}f(a),
\end{align*}
and set
\begin{equation}\label{eq:U}
 U_n=\bigcap_{j=1}^n\ker A_j\cap\bigcap_{j=1}^n\ker E_j.
\end{equation}
In other words, $f\in U_n$ if it vanishes whenever two adjacent labels
are equal, including $a_n=a_1$, and its sum over any one coordinate
is zero when all other coordinates are fixed. We will show that
$F_n$ represents the relative commutant $P_n$ and that $F_n\cap U_n$
represents the common kernel of the cap maps.

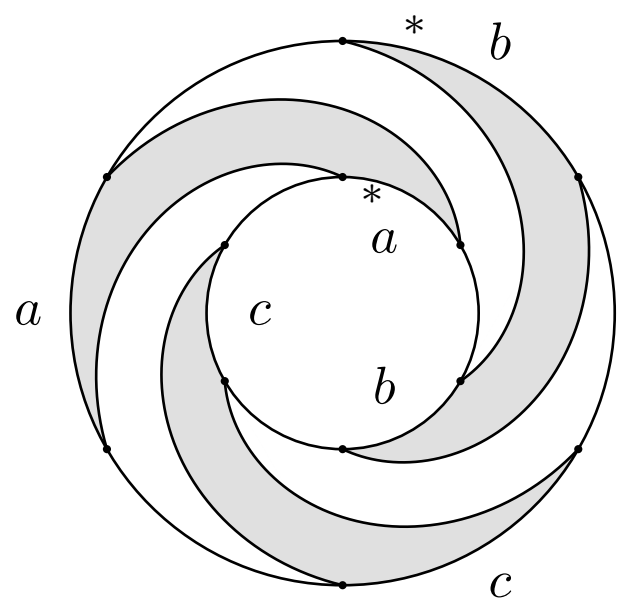
\begin{figure}[htbp]
\centering
\begin{tikzpicture}[x=1cm,y=1cm,>=stealth,
  every node/.style={font=\small},line width=.5pt]
% Each shaded strip joins the input interval to the output interval
% two boundary points counterclockwise from it.
\foreach \ang in {90,-30,-150}{
  \path[fill=black!12]
    plot[domain=0:1,samples=65,variable=\t]
      ({(.9+.9*\t)*cos(\ang+120*\t)},
       {(.9+.9*\t)*sin(\ang+120*\t)})
    arc[start angle=\ang+120,end angle=\ang+60,radius=1.8]
    plot[domain=1:0,samples=65,variable=\t]
      ({(.9+.9*\t)*cos(\ang-60+120*\t)},
       {(.9+.9*\t)*sin(\ang-60+120*\t)})
    arc[start angle=\ang-60,end angle=\ang,radius=.9] -- cycle;
}
\foreach \ang in {90,30,-30,-90,-150,150}{
  \draw plot[domain=0:1,samples=65,variable=\t]
      ({(.9+.9*\t)*cos(\ang+120*\t)},
       {(.9+.9*\t)*sin(\ang+120*\t)});
  \fill (\ang:.9) circle (.027);
  \fill (\ang:1.8) circle (.027);
}
\draw (0,0) circle (.9);
\draw (0,0) circle (1.8);
\foreach \ang/\lab in {60/a,-60/b,180/c}
  \node at (\ang:.55) {$\lab$};
\foreach \ang/\lab in {60/b,-60/c,180/a}
  \node at (\ang:2.08) {$\lab$};
\node[font=\scriptsize] at (76:.81) {$*$};
\node[font=\scriptsize] at (76:1.96) {$*$};
\node at (4.6,.42) {$e_{(a,b,c)}\ \longmapsto\ e_{(b,c,a)}$};
\node at (4.6,-.25) {$Z_{R_3}=\rho$};
\end{tikzpicture}
\caption{Two-click rotation in degree three. Gray regions are shaded,
and the stars mark the first shaded interval on the two boundaries.
Labels are read clockwise from these intervals. Following the shaded
strips from the inner boundary to the outer boundary gives the cyclic
permutation $(a,b,c)\mapsto(b,c,a)$.}
\label{fig:rotation}
\end{figure}

We first describe rotation in the spin planar algebra. An $n$-box
has $2n$ boundary points. We use the spin basis of
\cite[Definition~1.3]{M21}, reading the shaded intervals clockwise
from the marked shaded interval, as in Figure~\ref{fig:rotation}.
Note they are read counterclockwise starting after the marked 
interval in \cite{M21}.
Reversing the order interchanges $\rho$ and $\rho^{-1}$.
In our convention, $e_{(a_1,\ldots,a_n)}$ assigns the labels
$a_1,\ldots,a_n$ to the successive shaded intervals.
Rotation by two boundary points preserves the shading and moves each
label to the preceding shaded interval. This is called two-click
rotation. Figure~\ref{fig:rotation} shows this operation in degree three.

\begin{lemma}\label{lem:spinrotation}
Let $R_n$ be the two-click rotation tangle, oriented so that input
shaded interval $j+1$ is joined to output shaded interval $j$.
Its action in spin coordinates is
\begin{equation}\label{eq:spinrotation}
 (Z_{R_n})_{a,b}=\prod_{j=1}^n\delta_{a_j,b_{j+1}},
 \qquad Z_{R_n}=\rho,
\end{equation}
where $b_{n+1}=b_1$. Rotation in the opposite direction acts as $\rho^{-1}$.
\end{lemma}
\begin{proof}
Each shaded region of the rotation annulus joins one input interval
to one output interval. A spin state must be constant on this region.
Thus, for prescribed input labels $b$ and output labels $a$, a state
exists precisely when $a_j=b_{j+1}$ for all $j$, and in that case it
is unique.

In the spherical normalization of
\cite[Definition~1.3(iii), (iv), and (vi)]{M21}, the output contributes
a factor $d^{n/2}$, while each shaded region $r$ contributes
$d^{-\operatorname{Rot}(r)/2}$. The $n$ shaded regions are discs
whose positively oriented boundaries have rotation number one.
The coefficient of the unique state is therefore
\[
 d^{n/2}\prod_{j=1}^n d^{-1/2}=1.\qedhere
\]
\end{proof}

\begin{proposition}\label{prop:realization}
For every $n\geq1$, we have $T_n=C_nC_n^*$, and $T_n$ is a positive
contraction. There is a vector space isomorphism
$\mathcal J_n:P_n\to F_n$ satisfying
\begin{equation}\label{eq:rotationintertwiner}
 \mathcal J_n\operatorname{Rot}_{P_n}=\rho\mathcal J_n,
\end{equation}
where $\operatorname{Rot}_{P_n}$ is two-click planar rotation with the orientation of Lemma~\ref{lem:spinrotation}. Let $W_n\subset P_n$ be the common kernel of all elementary cap maps.
Under the above identification,
\[
 W_n=F_n\cap U_n\qquad(n\ge2).
\]
Also $W_1=0$.
\end{proposition}
\begin{proof}
We divide the proof into the description of the relative commutants
and the identification of their rotation and cap maps.

\emph{The relative commutants.}
In the finite dimensional tower of
\cite[Propositions~2.2--2.4]{N10}, write
\[
 \cP_{2r}=M_d(\C)^{\otimes(r+1)},\quad
 \cP_{2r+1}=M_d(\C)^{\otimes(r+1)}\otimes D_d,\quad
 \cP_{-1}=D_d.
\]
Let $U_i$ be the unitaries in that construction, with $U_0=U$.
Then
\[
 \cQ_i=U_i\cP_{i-1}U_i^*,\qquad \cA_i=D_d'\cap\cP_i,
 \qquad P_n\cong\cA_{n-1}\cap\cQ_{n-1}.
\]
Here $D_d$ is embedded in the first tensor factor. Since
$D_d'\cap M_d(\C)=D_d$, we have
\[
 \cA_{2r}=D_d\otimes M_d(\C)^{\otimes r},\qquad
 \cA_{2r+1}=D_d\otimes M_d(\C)^{\otimes r}\otimes D_d.
\]
Consequently,
\[
 \dim\cA_i=d^{i+1}=\dim\cP_{i-1}=\dim\cQ_i.
\]
Thus both $\cA_{n-1}$ and $\cQ_{n-1}$ admit orthonormal bases
indexed by $X^n$. We shall show that $T_n$ represents orthogonal
projection onto $\cQ_{n-1}$ followed by projection back onto
$\cA_{n-1}$. To compute its matrix, we first compute the inner
products between these bases.

The formulas for the unitaries give
\begin{align*}
 (U_{2r})_{x,y}
 &=d^{-(r+1)/2}\prod_{j=0}^rH_{x_jy_j}
                       \prod_{j=0}^{r-1}H_{x_{j+1}y_j},\\
 (U_{2r+1})_{x,y}
 &=\delta_{x_{r+1},y_{r+1}}d^{-(r+1)/2}
                       \prod_{j=0}^rH_{x_jy_j}H_{x_{j+1}y_j}.
\end{align*}
We use the unnormalized Hilbert--Schmidt inner product on each
ambient full matrix algebra. Let $E_{ab}$ denote a matrix unit.
At level $i=2r$, the matrix units of $\cA_i$ are indexed by words
of length $2r+1$ as follows:
\[
 (a_0,a_1,\ldots,a_r,b_r,\ldots,b_1)
 \longleftrightarrow
 E_{a_0a_0}\otimes E_{a_1b_1}\otimes\cdots\otimes E_{a_rb_r}.
\]
For $r\geq1$, index the corresponding basis of $\cQ_i$ by
$(c_0,\ldots,c_{r-1},s,d_{r-1},\ldots,d_0)$.
The inner product of a basis vector of $\cA_i$ with one of
$\cQ_i$ is the product of the corresponding two entries of $U_{2r}$.
By the formula above, the matrix of these inner products is $C_{2r+1}$.
At level $i=2r+1$, the two cyclic words are
\[
 (a_0,\ldots,a_r,t,b_r,\ldots,b_1),\qquad
 (c_0,\ldots,c_r,d_r,\ldots,d_0).
\]
The first word indexes
$E_{a_0a_0}\otimes E_{a_1b_1}\otimes\cdots\otimes E_{a_rb_r}\otimes E_{tt}$.
An orthonormal basis of the embedded algebra $\cP_{2r}$ is given
by $d^{-1/2}E_{c,d}\otimes I_d$, where
$E_{c,d}=E_{c_0d_0}\otimes\cdots\otimes E_{c_rd_r}$.
After conjugating by $U_{2r+1}$, its inner products with the above
basis of $\cA_i$ give $C_{2r+2}$. The factor $d^{-1/2}$ accounts
for the normalization in \eqref{eq:C}.

Let $\iota_n:\cH_n\to\cA_{n-1}$ denote these matrix-unit
coordinate maps. At the initial level,
$D_d\cap UD_dU^*=\C I_d$ and $C_1=J_d/d$, where $J_d$ is
the $d\times d$ matrix with every entry equal to one.

Put $A=\cA_{n-1}$ and $B=\cQ_{n-1}$, and let $E_A,E_B$ be
their orthogonal projections in the ambient Hilbert--Schmidt space.
Regard $\iota_n$ and the analogous coordinate map
$\eta_n:\cH_n\to B$ as isometries into this ambient space.
The preceding inner-product calculation says precisely that
$C_n=\iota_n^*\eta_n$. Since $E_B=\eta_n\eta_n^*$, we obtain
\[
 C_nC_n^*=\iota_n^*E_B\iota_n,\qquad
 \iota_nC_nC_n^*=E_AE_B\iota_n.
\]
In particular, $C_nC_n^*$ is a positive contraction. Moreover,
summing independently over the coordinates of $t\in X^n$ in
\eqref{eq:C} gives
\[
\begin{aligned}
 (C_nC_n^*)(a,b)
 &=d^{-n-1}\sum_{t\in X^n}\prod_{j=1}^n
 H_{a_jt_j}H_{a_{j+1}t_j}H_{b_jt_j}H_{b_{j+1}t_j}\\
 &=d^{-n-1}\prod_{j=1}^n
 \left(\sum_{t_j\in X}
 H_{a_jt_j}H_{b_jt_j}H_{a_{j+1}t_j}H_{b_{j+1}t_j}\right)\\
 &=d^{-n-1}\prod_{j=1}^n S(a_j,b_j,a_{j+1},b_{j+1})
 =T_n(a,b).
\end{aligned}
\]
This explains the formula \eqref{eq:T}.

To identify its fixed points, let $f\in\cH_n$ and put
$x=\iota_nf\in A$. Then
\[
 \langle(I-T_n)f,f\rangle
 =\|x\|^2-\|E_Bx\|^2
 =\|x-E_Bx\|^2.
\]
Thus $T_nf=f$ implies $x\in B$, and hence $x\in A\cap B$.
Conversely, if $x\in A\cap B$, both projections fix $x$, so
$T_nf=f$. Therefore
\[
 \iota_n\bigl(\ker(I-T_n)\bigr)
 =\cA_{n-1}\cap\cQ_{n-1}\cong P_n,
\]
which proves the required identification with $F_n$.

\emph{Rotation.}
To identify the planar operations on these spaces, we use the spin
construction of \cite[Example 2.8 and Propositions 2.8.1--2.8.2]{J99} and the construction of the flat subalgebra in \cite[Theorem 2.11.3 and Proposition 2.11.6]{J99}. These realize the Hadamard standard invariant as a flat shaded
planar subalgebra of the spin planar algebra. The identification with the commuting-square relative commutants is given in \cite[Proposition 2.11.12 and the following paragraph]{J99}; the corresponding finite dimensional description is given in \cite[Propositions 1.7--1.8 and Remark 1.9]{M21}. The identification is one of shaded planar algebras. In particular,
rotation and caps on the flat subalgebra are the restrictions of the
corresponding operations in the spin planar algebra.

Choose the spin component whose marked interval is shaded, and read
the labels around the boundary from that interval. Cutting the
boundary loop at its midpoint gives the matrix-unit coordinates
above: the first half supplies row indices and the second half
supplies column indices in reverse order. The initial label, and
also the midpoint label when $n$ is even, supplies a diagonal factor.
For example, in degree three we obtain
\[
 (a,b,c)\longleftrightarrow E_{aa}\otimes E_{bc}.
\]
The reading order is illustrated in Figure~\ref{fig:matrix-coordinates}.

\begin{figure}[htbp]
\centering
\begin{tikzpicture}[x=1cm,y=1cm,>=stealth,
  every node/.style={font=\small},line width=.5pt]
% A shaded collar records the boundary intervals, not interior strings.
\foreach \ang in {120,0,-120}{
  \path[fill=black!12]
    (\ang:1.24) arc[start angle=\ang,end angle=\ang-60,radius=1.24]
    -- (\ang-60:.82)
    arc[start angle=\ang-60,end angle=\ang,radius=.82] -- cycle;
}
\draw[densely dashed] (0,1.24)--(0,-1.24);
\draw (0,0) circle (1.24);
\foreach \ang in {120,60,0,-60,-120,180}
  \fill (\ang:1.24) circle (.03);
\node at (90:1.55) {$a$};
\node at (-30:1.02) {$b$};
\node at (210:1.02) {$c$};
\node[font=\scriptsize] at (109:1.37) {$*$};
\draw[->] (58:1.54) arc[start angle=58,end angle=-58,radius=1.54];
\draw[->] (122:1.54) arc[start angle=122,end angle=238,radius=1.54];
\node[right] at (1.68,0) {row};
\node[left] at (-1.68,0) {column};
\node at (0,-1.65) {$(a,b,c)$};
\draw[->] (2.65,0)--(3.45,0);
\node[anchor=west] at (3.85,.65) {row indices: $(a,b)$};
\node[anchor=west] at (3.85,.05) {column indices: $(a,c)$};
\node[anchor=west] at (3.85,-.75)
  {$E_{(a,b),(a,c)}=E_{aa}\otimes E_{bc}$};
\end{tikzpicture}
\caption{The matrix-unit coordinates in degree three. The dashed
line indicates the cut at the marked interval and the opposite
interval; it is not a string. Reading both halves from the marked
interval gives the row indices $(a,b)$ and the column indices $(a,c)$.
The column half is therefore read in the reverse of the cyclic order.
The common initial label $a$ gives the diagonal factor $E_{aa}$.}
\label{fig:matrix-coordinates}
\end{figure}
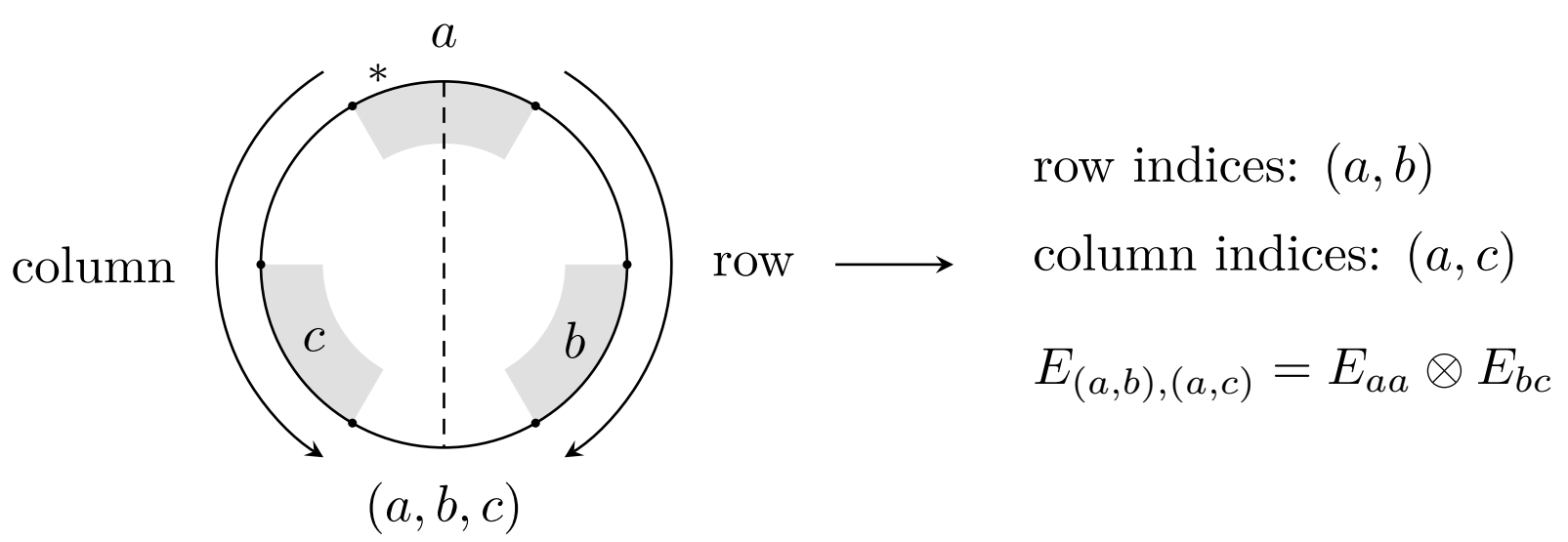

Since all spin labels have the same weight, the spherical spin basis
and these matrix units differ by a scalar depending only on the
degree. That scalar cancels in the formula for rotation.

Let $\Psi_n:P_n\to\cA_{n-1}\cap\cQ_{n-1}$ denote the
identification just described, and set
$\mathcal J_n=\iota_n^{-1}\Psi_n$. In the notation of
\cite[Remark~1.9 and Definition~1.10]{M21}, a flat vector $y$ in the
positive spin space has a partner $x$ in the negative spin space
satisfying $\psi_u(y)=1\otimes x$. In our orthonormal word bases this
is the equality $\eta_n y=\iota_n x$ in the intersection. Here
$\eta_n$ conjugates the indicated basis of $\cP_{n-2}$ by $U_{n-1}$, while
$\iota_n$ gives the matrix units on the $\cA_{n-1}$ side.
Applying $\iota_n^*$ and $\eta_n^*$ gives $x=C_n y$ and $y=C_n^*x$.
Thus $C_n$ maps $\ker(I-C_n^*C_n)$ unitarily onto $F_n$, with
inverse $C_n^*$ on $F_n$. By \eqref{eq:C}, simultaneous cyclic
permutation of the two words leaves $C_n(a,b)$ unchanged, so
$C_n\rho=\rho C_n$. Transporting the rotation of the flat vector $y$
to its partner $x$ by $C_n$, Lemma~\ref{lem:spinrotation} gives
\[
 \mathcal J_n\operatorname{Rot}_{P_n}\mathcal J_n^{-1}
 =Z_{R_n}|_{F_n}=\rho|_{F_n}.
\]
In degree three, the formula for the ambient rotation is thus
\[
E_{aa}\otimes E_{bc}\longmapsto E_{bb}\otimes E_{ca}.
\]
This proves \eqref{eq:rotationintertwiner}. The planar algebra
identification is essential here: the fact that $T_n$ commutes with
$\rho$ by itself only shows that $F_n$ is invariant under cyclic
permutation.

\begin{figure}[htbp]
\centering
\begin{tikzpicture}[x=1cm,y=1cm,
  every node/.style={font=\small},line width=.5pt]
\begin{scope}
  \path[fill=black!12]
    (1,0)..controls(1,1.55)and(2.6,1.55)..(2.6,0)--cycle;
  \draw (0,0)--(3.6,0);
  \draw (1,0)..controls(1,1.55)and(2.6,1.55)..(2.6,0);
  \fill (1,0) circle (.03); \fill (2.6,0) circle (.03);
  \node at (1.8,.48) {$x$};
  \node at (1.8,1.6) {$L_j^*$};
  \node[below] at (1.8,-.04) {input boundary};
  \node at (1.8,-1.15) {$\displaystyle\sum_{x\in X}f(\ldots,x,\ldots)$};
\end{scope}
\begin{scope}[shift={(6,0)}]
  \fill[black!12] (0,0) rectangle (3.6,1.35);
  \path[fill=white]
    (1,0)..controls(1,1.55)and(2.6,1.55)..(2.6,0)--cycle;
  \draw (0,0)--(3.6,0);
  \draw (1,0)..controls(1,1.55)and(2.6,1.55)..(2.6,0);
  \fill (1,0) circle (.03); \fill (2.6,0) circle (.03);
  \node at (.45,.43) {$a$}; \node at (3.15,.43) {$b$};
  \node at (1.8,1.6) {$D_j^*$};
  \node[below] at (1.8,-.04) {input boundary};
  \node at (1.8,-1.15) {$f(\ldots,a,a,\ldots)$};
\end{scope}
\end{tikzpicture}
\caption{The two elementary cap operations, shown locally along the
input boundary. On the left, the cap closes a shaded region, and its
label $x$ is summed. On the right, the two shaded intervals belong to
one connected region, so their labels satisfy $a=b$. The displayed
maps are $L_j^*$ and $D_j^*$ for the counting inner products; the
planar cap maps are nonzero scalar multiples of these maps. All other
labels are unchanged.}
\label{fig:caps}
\end{figure}
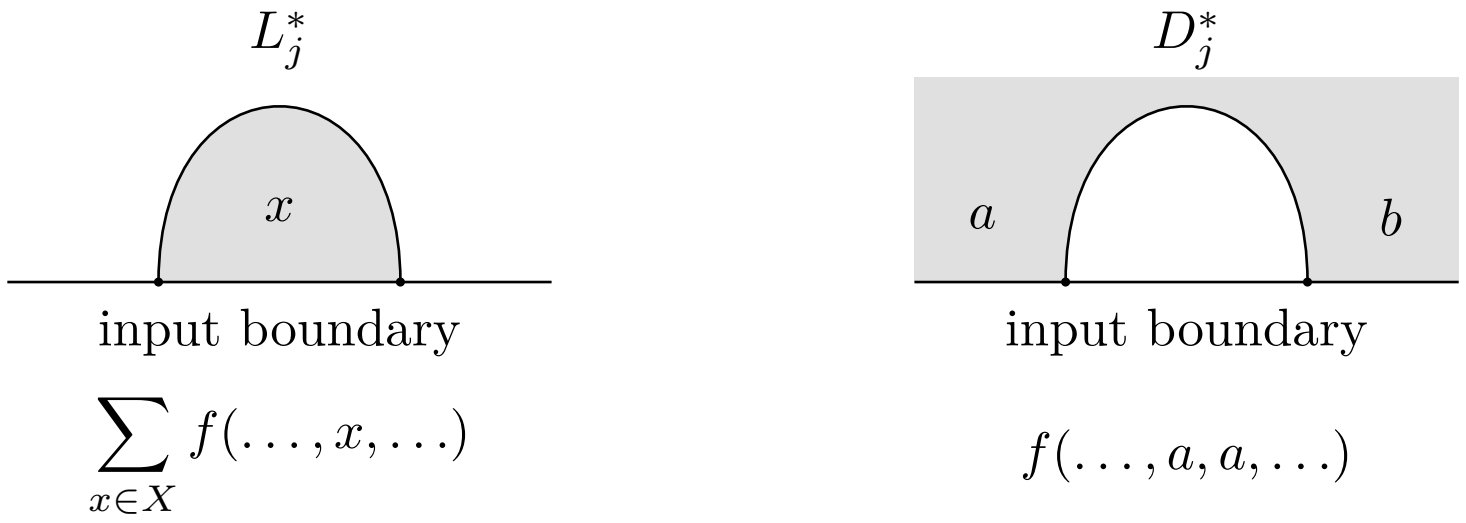

\emph{Caps.}
In spin coordinates a cup either introduces a free label or repeats
an existing label. The corresponding cap operations are shown in
Figure~\ref{fig:caps}. Let $L_j:\cH_{n-1}\to\cH_n$ introduce a free
$j$th coordinate, so that $(L_jf)(a)$ is obtained by omitting $a_j$
from the arguments of $f$. Let $D_j$ repeat a label in adjacent
positions $j,j+1$, with value zero when the two labels differ.
Then
\[
 A_j=d^{-1}L_jL_j^*,\qquad E_j=D_jD_j^*,\qquad
 L_j^*L_j=dI,\quad D_j^*D_j=I.
\]
In each degree, the trace inner product is a positive scalar
multiple of the counting inner product. The same holds for the
spherical inner product by \cite[Theorem~4.2.1(iv)]{J99}.
The cap maps are therefore nonzero scalar multiples of $L_j^*$
and $D_j^*$. The equation $L_j^*f=0$ says that the sum over the
$j$th coordinate is zero. The equation $D_j^*f=0$ says that $f$
vanishes when the corresponding adjacent labels are equal.
By \eqref{eq:intertwineC} and its adjoints, $C_n$ and $C_n^*$
each carry $U_n$ into itself. Thus the flat vector $y$ lies in the
common cap kernel precisely when its partner $x=C_ny$ lies in $U_n$.
This proves $W_n=F_n\cap U_n$ for $n\geq2$.
Finally, $F_1$ consists of constant functions, and the cap does not
vanish on a nonzero constant. Thus $W_1=0$.
\end{proof}

The preceding description can also be checked directly from the
matrices. Row and column orthogonality give
\begin{equation}\label{eq:intertwineC}
 C_nA_j=E_jC_n,\qquad C_nE_j=A_{j+1}C_n.
\end{equation}
Together with their adjoints, these identities show that $T_n$
commutes with each $A_j$ and $E_j$. Also, row and column
orthogonality give
\begin{equation}\label{eq:contractions}
 S(a,b,c,b)=d\delta_{ac},\qquad
 \sum_xS(a,b,c,x)S(c,x,e,f)=dS(a,b,e,f)
\end{equation}
and consequently
\[
T_nD_j=D_jT_{n-1},\qquad T_nL_j=L_jT_{n-1}.
\]
For the second identity in \eqref{eq:contractions}, expand the two
sums defining $S$ and sum first over $x$. It follows that cups and
caps preserve the spaces $F_n$ in their respective degrees.
Equation~\eqref{eq:T} shows that $T_n$ commutes with $\rho$, so
both $F_n$ and $W_n$ are invariant under rotation.

\FloatBarrier

\section{A congruence and its consequence for rotation}\label{sec:parity}

Throughout this section we assume that $d=4m>4$, where $m$ is odd.
We first record a congruence satisfied by the sums in
\eqref{eq:profile}. We then use it to prove that rotation is the
identity on $W_n$ whenever $n$ is odd.

The following four-row congruence is classical
\cite[p.~425]{W89}. Montgomery obtains the identity
$S(a,b,c,e)=8k-d$, where $k$ counts the columns specified below,
in the proof of \cite[Proposition~3.2.9]{M22}.

\begin{lemma}\label{lem:parity}
If $a,b,c,e$ are four distinct row labels, then
\[
 S(a,b,c,e)\equiv-d\pmod8.
\]
In particular, $S(a,b,c,e)/4$ is an odd integer.
\end{lemma}
\begin{proof}
Multiply column $t$ by $H_{at}$. This preserves row inner products
and products of four row entries, since $H_{at}^2=H_{at}^4=1$.
Row $a$ now consists of ones, and the other three rows have entries
\[
 x_t=H_{at}H_{bt},\qquad
 y_t=H_{at}H_{ct},\qquad
 z_t=H_{at}H_{et}.
\]
Orthogonality to row $a$ and pairwise orthogonality give, respectively,
\[
 \sum_t x_t=\sum_t y_t=\sum_t z_t=0,
 \qquad
 \sum_t x_ty_t=\sum_t x_tz_t=\sum_t y_tz_t=0.
\]
Moreover, $H_{at}^3=H_{at}$, so
\[
 \sum_t x_ty_tz_t
 =\sum_t H_{at}H_{bt}H_{ct}H_{et}
 =S(a,b,c,e).
\]
Expanding the product and using these identities, we obtain
\[
\begin{aligned}
 \sum_t(1+x_t)(1+y_t)(1+z_t)
 &=d+\sum_t(x_t+y_t+z_t)\\
 &\quad+\sum_t(x_ty_t+x_tz_t+y_tz_t)
       +\sum_t x_ty_tz_t\\
 &=d+S(a,b,c,e).
\end{aligned}
\]
Each summand on the left is eight when $x_t=y_t=z_t=1$ and zero
otherwise. Thus $d+S(a,b,c,e)$ is eight times the number of such
columns, proving the congruence. Since $d\equiv4\pmod8$, it follows
that $S(a,b,c,e)\equiv4\pmod8$, so $S(a,b,c,e)/4$ is odd.
\end{proof}

For two-element subsets of $X$, define
\[
 \kappa(\{a,c\},\{b,e\})=\frac14S(a,b,c,e).
\]
If the two pairs are equal, $\kappa=m$. If they have exactly one
common element, orthogonality gives $\kappa=0$. If they are
disjoint, Lemma~\ref{lem:parity} shows that $\kappa$ is an odd
integer. These three cases can be written as the single congruence
\begin{equation}\label{eq:localmod2}
 \kappa(\{a,c\},\{b,e\})
 \equiv1+\delta_{ab}+\delta_{ae}+\delta_{cb}+\delta_{ce}\pmod2
 \qquad(a\ne c,\ b\ne e).
\end{equation}

\begin{lemma}\label{lem:firstdegrees}
We have $\dim F_1=1$ and $\dim F_2=2$.
\end{lemma}
\begin{proof}
The first assertion follows from $T_1=J_d/d$.
The second is the two-dimensionality result for the second
relative commutant in \cite[Proposition~3.2.9]{M22}. In the present coordinates,
\[
 T_2((a,c),(b,e))=d^{-3}S(a,b,c,e)^2.
\]
The subspace spanned by $e_{(a,a)}$, $a\in X$, is invariant under
$T_2$, and the restriction to it is $J_d/d$. The complementary
subspace, spanned by $e_{(a,c)}$ with $a\ne c$, is also invariant.
Let $B$ be the restriction of $T_2$ to this subspace.
It is symmetric and has nonnegative entries. For $a\ne c$,
column orthogonality gives
\[
 \sum_{b,e}S(a,b,c,e)^2=d^3.
\]
Thus each row sum of $B$ is one. Consider the graph whose vertices
are ordered pairs $(a,c)$ with $a\ne c$, with an edge whenever
the corresponding entry of $B$ is positive. The vertices $(a,c)$
and $(c,a)$ are joined. Pairs with disjoint underlying sets are
joined by Lemma~\ref{lem:parity}. If two distinct underlying pairs
meet, there is a two-element set disjoint from their union, since
$d>4$. The graph is therefore connected.

If $Bf=f$, symmetry and the row-sum condition imply
\[
 0=\langle(I-B)f,f\rangle
   =\frac12\sum_{u,v}B_{uv}|f(u)-f(v)|^2.
\]
It follows that $f$ is constant on this connected graph.
Each of the two invariant subspaces therefore contributes exactly
one dimension to $F_2$.
\end{proof}

For $n\geq2$, let
\[
\Omega_n=\{a\in X^n:a_j\ne a_{j+1}\text{ for every }j\}.
\]
Every vector in $U_n$ is supported on $\Omega_n$.
Define an integer matrix indexed by this set by
\[
 K_n(a,b)=\prod_{j=1}^n
       \kappa(\{a_j,a_{j+1}\},\{b_j,b_{j+1}\}).
\]
Since $S=4\kappa$ and $d=4m$, the submatrix of $T_n$ indexed
by $\Omega_n$ is $K_n/(dm^n)$. Thus, restricting a vector in
$W_n=F_n\cap U_n$ to $\Omega_n$, we have
\begin{equation}\label{eq:integerfixed}
 K_nf=dm^nf\qquad(f\in W_n).
\end{equation}

\begin{lemma}\label{lem:globalmod2}
For $n\ge2$ and every integer-valued $f\in U_n$,
\[
 K_nf\equiv(\rho+\rho^{-1})f\pmod2
\]
on $\Omega_n$.
\end{lemma}
\begin{proof}
Fix $a\in\Omega_n$. Since $f$ vanishes outside $\Omega_n$,
we may sum over all $b\in X^n$ when expanding
\eqref{eq:localmod2} in the formula for $(K_nf)(a)$.
Any term that does not involve a particular input coordinate $b_j$
sums to zero, because the sum of $f$ over that coordinate is zero.
Each summand selected from a factor involves at most one input coordinate.
Consequently, a term that remains must select one Kronecker delta
from each factor and use each input coordinate exactly once.

For such a term, put $\epsilon_j=0$ if the $j$th factor uses $b_j$,
and $\epsilon_j=1$ if it uses $b_{j+1}$. The number of occurrences
of $b_j$ is
\[
1-\epsilon_j+\epsilon_{j-1}.
\]
It is equal to one for every $j$, so all the $\epsilon_j$ are equal.
There are therefore only two possibilities: every factor uses its
first input coordinate, or every factor uses its second.

In either case, each selected delta equates that input coordinate
with $a_j$ or $a_{j+1}$. If these choices are not constant around
the cycle, there is a transition from the second choice to the first.
Two adjacent input labels are then equal, and the corresponding
value of $f$ is zero. The four constant choices give
$f(a),f(a),(\rho f)(a),(\rho^{-1}f)(a)$.
The first two cancel modulo two, proving the assertion.
\end{proof}

\begin{proposition}\label{prop:oddrotation}
If $n\geq3$ is odd, then $\rho|_{W_n}=I$.
\end{proposition}
\begin{proof}
All equations defining $W_n$ have rational coefficients. If
$W_{n,\Q}=W_n\cap\Q^{X^n}$, then
$W_n=W_{n,\Q}\otimes_{\Q}\C$. Consider the operator
\[
R=I+\rho+\cdots+\rho^{n-1}.
\]
Suppose that $W_{n,\Q}\cap\ker R$ contains a nonzero vector.
Clearing denominators and dividing by the greatest common divisor
of its coordinates gives an integer-valued vector $g$ in this
intersection with at least one odd coordinate.

By \eqref{eq:integerfixed}, $K_ng=dm^ng$. The right-hand side
is zero modulo two. Lemma~\ref{lem:globalmod2} therefore gives
$(\rho+\rho^{-1})g\equiv0\pmod2$, and hence
\[
 \rho^2g\equiv g\pmod2.
\]
Since $n$ is odd and $\rho^n=I$, we have
$\rho=(\rho^2)^{(n+1)/2}$. It follows that
$\rho g\equiv g\pmod2$. But then
\[
 0=Rg\equiv ng\equiv g\pmod2,
\]
contrary to the choice of $g$. Thus $W_{n,\Q}\cap\ker R=0$.
Since $W_{n,\Q}$ is invariant under $\rho$ and
$R(\rho-I)=\rho^n-I=0$, the image
$(\rho-I)W_{n,\Q}$ lies in this zero intersection.
Therefore $\rho=I$ on $W_{n,\Q}$, and then on $W_n$.
\end{proof}

\section{The trace of rotation}\label{sec:trace}

We next bound the trace of rotation on $P_n$, independently of $n$.
By Proposition~\ref{prop:realization}, this is the trace of $\rho$
on $F_n$. The estimate holds for every real Hadamard matrix;
no congruence assumption on its order is needed.
Throughout this section, $\Tr$ denotes the ordinary matrix trace.

\begin{proposition}\label{prop:trace}
For every $n\ge1$,
\[
 |\Tr(\rho|_{F_n})|\le d.
\]
\end{proposition}
\begin{proof}
Set
\[
 \xi_{ab}=d^{-1/2}(H_{ay}H_{by})_{y\in X},\qquad
 Q_{ab}=|\xi_{ab}\rangle\langle\xi_{ab}|,\qquad
 \tau_d=d^{-1}\Tr.
\]
Here $Q_{ab}$ is the orthogonal projection onto $\C\xi_{ab}$.
For either fixed index, the vectors $\xi_{ab}$ form an orthonormal
basis. Thus the corresponding projections are mutually orthogonal
and sum to $I_d$. It follows that the block matrix
$V=(Q_{ab})\in M_d(\C)\otimes M_d(\C)$ is unitary. Indeed,
\[
(VV^*)_{ac}=\sum_bQ_{ab}Q_{cb}=\delta_{ac}I_d,
\]
and the same calculation gives $V^*V=I$.

Taking the trace of a product of rank-one projections yields
\begin{equation}\label{eq:Qtrace}
 T_n(a,b)=\tau_d(Q_{a_1b_1}\cdots Q_{a_nb_n}).
\end{equation}
For $r\geq1$, number the tensor factors by $0,1,\ldots,r$.
Let $V_{0,j}$ act as $V$ on factors $0,j$ and as the identity on
the other factors. The operator
\[
 V^{[r]}=V_{0,1}V_{0,2}\cdots V_{0,r}.
\]
is unitary. Viewed as a block matrix with respect to factor zero,
it has entries
\[
 V^{[r]}_{ab}=\sum_{c_1,\ldots,c_{r-1}}
 Q_{ac_1}\otimes Q_{c_1c_2}\otimes\cdots\otimes Q_{c_{r-1}b}.
\]
Put $\tau_r=\tau_d^{\otimes r}$. Applying \eqref{eq:Qtrace} to
each factor in the matrix product $T_n^r$, we obtain
\begin{equation}\label{eq:powers}
 (T_n^r)(a,b)=\tau_r(V^{[r]}_{a_1b_1}\cdots V^{[r]}_{a_nb_n}).
\end{equation}
Indeed, the sums over $c_1,\ldots,c_{r-1}$ in each block entry
are precisely the sums over intermediate words in $T_n^r$.
Formulas \eqref{eq:Qtrace} and \eqref{eq:powers} are the
transfer-matrix formulas of Banica \cite[Proposition~3.2]{B14}
for the projections $Q_{ab}$.
Taking the trace after applying $\rho$, we have
\begin{align*}
 \Tr(\rho T_n^r)
 &=\sum_{a_1,\ldots,a_n}
   \tau_r(V^{[r]}_{a_na_1}V^{[r]}_{a_1a_2}\cdots V^{[r]}_{a_{n-1}a_n})\\
 &=\tau_r\bigl(\Tr_d((V^{[r]})^n)\bigr).
\end{align*}
Here $\Tr_d$ is the ordinary trace on factor zero. After division
by $d$, the last expression is the normalized trace of the unitary
$(V^{[r]})^n$. Hence $|\Tr(\rho T_n^r)|\leq d$.

Since $T_n$ is a positive contraction, its powers converge to the
orthogonal projection onto $F_n$. This projection commutes with
$\rho$. Letting $r\to\infty$ proves the desired estimate.
\end{proof}

\section{Annular multiplicities}\label{sec:annular}
In this section, let $P$ be a connected irreducible spherical
$C^*$-planar algebra with modulus $\delta>2$. Let $W_n$ be the
common kernel of the elementary cap maps on $P_n$, and put
\[
b_n=\dim P_n,\qquad a_n=\dim W_n.
\]
The numbers $a_n$ are the annular multiplicities. They count the
lowest-weight vectors of degree $n$ in the annular decomposition.
We recall the form of Jones's theorem that we will use.

\begin{theorem}[Jones's annular decomposition]\label{thm:annular}
For $k\geq1$ and $\omega^k=1$, there is an irreducible annular
Temperley--Lieb module $V^{k,\omega}$ whose lowest-weight space
is one-dimensional, lies in degree $k$, and has rotation eigenvalue
$\omega$. Its dimensions are
\[
 \dim V_n^{k,\omega}=\binom{2n}{n-k}\quad(n\ge k),\qquad
 V_n^{k,\omega}=0\quad(n<k).
\]
As an annular Temperley--Lieb module, $P$ has an orthogonal decomposition
\[
 P\cong V^{\mathrm{vac}}\oplus
       \bigoplus_{k\ge1}\ \bigoplus_{\omega^k=1}
       (V^{k,\omega})^{\oplus m_{k,\omega}},
 \qquad m_{k,\omega}=\dim\bigl(W_k\cap\ker(\rho-\omega I)\bigr).
\]
Here $V^{\mathrm{vac}}$ is the module generated by the zero-box
spaces, and $V_n^{\mathrm{vac}}=\TL_n(\delta)$.
\end{theorem}
We use this as a quoted result of Jones; see
\cite[Corollary~3.9, Definition~5.2, Proposition~5.3,
Theorem~5.7, and Corollary~5.8]{J01}, together with the
orthogonal decomposition in \cite[Remark~3.6]{J01}.
Connectedness implies that the vacuum module contains both
zero-box spaces, so no other summand of weight zero occurs.
Irreducibility gives $a_1=0$. Taking dimensions, we obtain
\begin{equation}\label{eq:annulardim}
 b_n=\Cat_n+\sum_{k=2}^n\binom{2n}{n-k}a_k.
\end{equation}
The same identity follows from \cite[Theorem~6.3]{J01}.

We will also use the basis described in \cite[Definition~5.2]{J01}.
In degree $n$, a module of weight $k>0$ has a basis of annular
link diagrams with $2k$ strings joining the inner and outer
boundaries. Diagrams that differ by inner rotation are identified
with the scalar factor prescribed by $\omega$. Linear independence
follows from the positivity theorem cited above. Outer rotation
permutes these basis classes, possibly multiplying them by scalars.
Thus only classes fixed by outer rotation can contribute to its trace.

\begin{lemma}\label{lem:oddcharacter}
For every odd $n\ge3$,
\[
 \Tr(\rho|_{P_n})=2+\Tr(\rho|_{W_n}).
\]
\end{lemma}
\begin{proof}
Consider first a summand of weight $k$ with $0<k<n$.
A basis diagram determines the set of $2k$ outer endpoints of its
through strings. This set is unchanged by inner rotation. Therefore,
if a basis class is fixed by two-click outer rotation, its endpoint
set must also be fixed.

On the $2n$ boundary points, two-click rotation has two orbits of
length $n$. Any invariant subset has cardinality $0,n$, or $2n$.
Since $n$ is odd, none of these can equal the even number $2k$
with $0<2k<2n$. Thus every summand of positive weight less than
$n$ contributes zero to the trace.

On the vacuum module, rotation permutes the noncrossing pairings.
Each such pairing contains an adjacent pair. If the pairing is
invariant under rotation, it contains the full orbit of this pair.
The orbit has length $n$: a stabilizer of an unordered pair has
order at most two, whereas the rotation group has odd order.
The orbit therefore accounts for all $n$ pairs. There are exactly
two such pairings, obtained by pairing alternating adjacent boundary
points. The vacuum contribution to the trace is two.

Finally, the degree-$n$ spaces of all summands of weight $n$ form
$W_n$. Their contribution is $\Tr(\rho|_{W_n})$, as required.
\end{proof}

We now return to the subfactor associated to a real Hadamard matrix
of order $d>4$ with $d\equiv4\pmod8$.
For odd $n\geq3$, Proposition~\ref{prop:oddrotation} gives
$\Tr(\rho|_{W_n})=a_n$. The preceding lemma and
Proposition~\ref{prop:trace} imply $2+a_n\leq d$. Hence
\begin{equation}\label{eq:bounded}
 0\le a_n\le d-2\qquad(n\ge3\text{ odd}).
\end{equation}

\section{Bounded annular multiplicities and graph norm}\label{sec:norm}
The estimate \eqref{eq:bounded} has a consequence for the principal
graph that does not depend on the Hadamard construction. We prove
this consequence first, using the dimensions of relative commutants
and the spectral measure of the principal graph at its distinguished
vertex.

\begin{proposition}\label{prop:norm}
Let $P$ be the planar algebra of an irreducible finite-index subfactor
with modulus $\delta>2$. For each $n\geq 1$, let $a_n$ be the $n$th
annular multiplicity of $P$, that is, the dimension of the subspace
of $P_{n,+}$ annihilated by every elementary cap map. If there is a
constant $C$ such that $a_n\leq C$ for infinitely many $n$, then the
principal graph of the subfactor is $A_\infty$, rooted at its endpoint.
\end{proposition}

\begin{proof}
Let $\Gamma$ be the principal graph, with distinguished vertex $*$
and adjacency operator $A$. Its positive dimension function
$v\mapsto d_v$ satisfies $\sum_w A_{vw}d_w=\delta d_v$.
The weighted Schur test gives $\|A\|\leq\delta$.
Let $\mu$ be the spectral measure of $A^2$ at the unit vector
$e_*$. It is a probability measure on $[0,\delta^2]$, and the
path-algebra description in \cite[Theorem~3.1.3]{J99} gives
\[
 b_j=\langle A^{2j}e_*,e_*\rangle
     =\int x^j\,d\mu(x).
\]

Consider the polynomials defined by $q_0(x)=2$, $q_1(x)=x-2$ and
\[
 q_{n+1}=(x-2)q_n-q_{n-1}.
\]
These are rescaled Chebyshev polynomials:
$q_n(2+2\cos\theta)=2\cos(n\theta)$.
Substituting $x=2+z+z^{-1}$ and expanding gives
\begin{equation}\label{eq:binomial}
 x^n=\binom{2n}{n}+\sum_{k=1}^n\binom{2n}{n-k}q_k(x).
\end{equation}
Since $b_0=b_1=1$, we have $\int q_1\,d\mu=-1$.
Integrating \eqref{eq:binomial}, the constant and $k=1$ terms combine to
$\binom{2n}{n}-\binom{2n}{n-1}=\Cat_n$.
Comparison with \eqref{eq:annulardim}, successively in degrees
$2,3,\ldots,n$, gives
\begin{equation}\label{eq:qmoment}
 a_n=\int q_n(x)\,d\mu(x)\qquad(n\geq2).
\end{equation}

The Chebyshev formula shows that $q_n(x)\geq-2$ for
$0\leq x\leq4$. For $x>4$, we have
\[
 q_n(x)=\lambda(x)^n+\lambda(x)^{-n},\qquad
 \lambda(x)=\frac{x-2+\sqrt{x(x-4)}}2>1.
\]
If $\mu((4,\delta^2])>0$, choose $c>4$ such that
$\varepsilon=\mu([c,\delta^2])>0$.
Since $\lambda$ is increasing and $q_n(x)\geq-2$ on
$[0,\delta^2]$, equation~\eqref{eq:qmoment} implies
\[
 a_n\geq\varepsilon\lambda(c)^n-2.
\]
This contradicts boundedness along an unbounded sequence.
Hence $\mu$ is supported on $[0,4]$, and therefore
\[
 b_j=\int x^j\,d\mu(x)\leq4^j\qquad(j\geq1).
\]

Finally, let $G$ be a finite connected subgraph containing $*$,
and let $v$ be a positive unit Perron eigenvector of its adjacency
matrix $A_G$. Counting closed walks and using the spectral
decomposition of $A_G$, we obtain
\[
 |v_*|^2\|A_G\|^{2j}\le(A_G^{2j})_{*,*}\le b_j\le4^j.
\]
Taking $2j$th roots and letting $j\to\infty$ gives
$\|A_G\|\leq2$. Since $\Gamma$ is connected, such subgraphs
exhaust $\Gamma$, and $\|\Gamma\|=\|A\|\leq2$.

If $\Gamma$ were finite, its positive dimension vector would be an
eigenvector of $A$ with eigenvalue $\delta>2$, contradicting the norm
bound. Thus $\Gamma$ is infinite. Irreducibility gives a root of
degree one, so the classification in \cite[Section~2.10]{P23} leaves
only $A_\infty$ and $D_\infty$. Both graphs have an $A_\infty$ tail.
Since the subfactor is extremal by irreducibility and has index
$\delta^2>4$, Popa's theorem on principal graphs with an $A_\infty$
tail \cite{P95} gives $\Gamma=A_\infty$, rooted at its endpoint.
\end{proof}

\begin{proof}[Proof of Theorem~\ref{thm:main}]
For every odd $n\geq3$, Proposition~\ref{prop:oddrotation} says
that rotation is the identity on $W_n$. Combining the trace estimate
of Proposition~\ref{prop:trace} with Lemma~\ref{lem:oddcharacter}
gives $0\leq a_n\leq d-2$, as in \eqref{eq:bounded}.
Proposition~\ref{prop:norm} therefore gives $\Gamma=A_\infty$.

Closed walks of length $2n$ starting at the endpoint of $A_\infty$
are counted by $\Cat_n$. Hence $\dim P_n=\Cat_n$ for every $n$.
For $\delta>2$, the Temperley--Lieb Markov form is positive definite;
see \cite[Section~4]{J01}. Thus the canonical Temperley--Lieb
subspace of $P_n$ already has dimension $\Cat_n$, and
\[
P_{n,+}=\TL_{n,+}(\delta).
\]
For $n\geq1$, one-click rotation is an invertible map from
$P_{n,+}$ to $P_{n,-}$ carrying Temperley--Lieb diagrams to
Temperley--Lieb diagrams. It follows that
$P_{n,-}=\TL_{n,-}(\delta)$ as well. In degree zero the equality
holds by connectedness. This identifies the entire shaded standard
invariant, and proves the assertions about both principal graphs
and the dimensions of the relative commutants.
\end{proof}

\Needspace{12\baselineskip}
\begin{corollary}\label{cor:rank}
Let $H$ be a real Hadamard matrix of order $d>4$ with
$d\equiv4\pmod8$, and let $T_n$ be the matrices defined in
\eqref{eq:T}. Then
\[
 \dim\ker(I-T_n)=\Cat_n,\qquad
 \operatorname{rank}(I-T_n)=d^n-\Cat_n\qquad(n\ge1).
\]
Moreover, $W_n=0$ for every $n\geq1$.
\end{corollary}
\begin{proof}
The dimension formula follows from Proposition~\ref{prop:realization}
and Theorem~\ref{thm:main}; the rank formula is then immediate.
Since $b_n=\Cat_n$, the nonnegative terms in
\eqref{eq:annulardim} must all vanish. Together with $W_1=0$,
this gives the last assertion.
\end{proof}

\end{document}